\documentclass[12pt,a4paper]{amsart}

\newcommand{\N}{\rm{I\!N}}        
\newcommand{\R}{\rm{I\!R}}

\newtheorem{theorem}{Theorem}

\newtheorem{coro}[theorem]{Corollary}
\newtheorem{lemma}[theorem]{Lemma}
\newtheorem{prop}[theorem]{Proposition}

\newcommand{\eps}{\varepsilon}

\newcommand{\norm}[1]{\|#1\|}

\newcommand{\betr}[1]{| #1  |}
\newcommand{\eing}[1]{_{|{#1}}}
\newcommand{\ebew}{\hfill{\rule{1.2ex}{1.2ex}}}
\newcommand{\bgl}{\begin{eqnarray}}
\newcommand{\bglst}{\begin{eqnarray*}}
\newcommand{\egl}{\end{eqnarray}}
\newcommand{\eglst}{\end{eqnarray*}}

\newcommand{\wst}{$w^{*}$}
\newcommand{\mdE}{;\,}
\newcommand{\ball}{{\mathcal B}}
\newcommand{\be}{\begin{enumerate}}
\newcommand{\ee}{\end{enumerate}}
\renewcommand{\limsup}{\overline{\lim}\;}

\newcommand{\wstAb}[1]{\overline{#1}^{*}}

\newcommand{\acoinfty}{\mbox{aco}_{\infty}}

\newcommand{\fkt}[5]{\begin{array}{rrcl}
                          #1: & #2 &\rightarrow&#3\\
                     \mbox{ } & #4 &\mapsto    &#5 \end{array}}
\newcommand{\Real}{\mbox{Re\,}}

\newcounter{Bem}
\newcounter{Fra}
\begin{document}
\title{A conjecture of Godefroy concerning James' theorem}
\author{Hermann Pfitzner}
\address{Universit\'e d'Orl\'eans\\
BP 6759\\
F-45067 Orl\'eans Cedex 2\\France}
\email{hermann.pfitzner@univ-orleans.fr}
\begin{abstract}
In this note we look at the interdependences between James' theorem and the boundary problem.
To do so we show a variant of James' sup-theorem for $C(K)$-spaces conjectured by Godefroy:
in order to know that a bounded weakly closed subset of a $C(K)$-space is weakly compact it is enough to test
the sup-attainment only for measures with countable support.
\end{abstract}
\maketitle\noindent
\noindent
{\bf \S 1. Introduction.}
This note is an afterthought to the solution of the boundary problem, proved in \cite{Pf-boundary}
and cited in Theorem \ref{boundary} below, and hopefully sheds some light
onto the interdependence between James's theorem \cite{James-sup-Th} and the boundary theorem. Both theorems provide a criterion for weak compactness.
Our idea is to show that, roughly speaking, the two criterions are, in some sense, equivalent and dual to each other.

Less roughly speaking,
let us consider an easy special case. (For notation and definitions see below.)
Let $X$ be a Banach space with unit ball $\ball(X)$ and dual $X^*$.
Suppose we know that each bounded linear functional on $X$ attains its norm on $\ball(X)$.
Then, on the one hand, James' theorem says that $\ball(X)$ is weakly compact and the boundary theorem (which we apply to the dual $X^*$)
says that the dual unit ball $\ball(X^*)$ is weakly compact (because by assumption $\ball(X)$ is a boundary for $X^*$ and because by
Banach-Alaoglu  $\ball(X^*)$ is compact in the topology of  pointwise convergence on $\ball(X)$).
On the other hand, it is well-known that $\ball(X)$ is weakly compact iff $\ball(X^*)$ is; this is due to Grothendieck's
classical double limit criterion which displays a kind of symmetry between weak compactness in $X$ and its dual
and which we take, in our context, as a generic theorem which includes related theorems like the one of Eberlein-\v{S}mulyan or of Gantmacher
-- or like Lemma \ref{Groth} below.
To resume, in the present situation James' theorem and the boundary theorem have equivalent conclusions namely
the weak compactness of $\ball(X)$ and, respectively, of $\ball(X^*)$.
Note, however, that the symmetry between both theorems is imperfect because in the example the boundary theorem yields James' theorem
but in the converse sense only the boundary theorem for the dual $X^*$ is recovered by James' theorem.
In order to study to which extent James' theorem implies the boundary theorem and vice versa, in the general case
(i.e.\ when the unit balls above are replaced by arbitrary bounded subsets),
it is necessary to find variants of both theorems.
We suggest two ways to do so.
First, we use the main result of \cite{Pf-boundary} (which is more general than the solution of the boundary problem)
in order to show a variant of James' theorem, Theorem \ref{James2.0} below,
which goes back to a conjecture of Godefroy \cite[p.~ 46]{G-Murcia}.
Godefroy's main interest in his conjecture was the possibility to
derive from it a positive solution to the boundary problem.
See \eqref{gl02} below for a summary of these interdependences.

The shortcoming of the just described implications in \eqref{gl02} is that they cannot be reversed entirely, for example
the stronger form of James' theorem is enough to yield the boundary theorem but cannot be recoverd entirely by the latter.
(This lack of symmetry has already been observed in the example above.)
In order to remedy this, a second way to present the interdependence of both theorems,
resumed in line \eqref{gl04} below, is suggested: Here the stronger form of James' theorem
is strengthend still a bit more (cf.\ Prop.\ \ref{James2.1}) so to become equivalent to the stronger form of the boundary theorem.
The inconvenience of this second variant \eqref{gl04} is that the statements of the involved theorems are less digest than
their 'classical' counterparts
because they are more technical but this is outweighed by two convincingly short proofs.\smallskip

Put differently and succintly, in this note only the proof of Godefroy's conjecture is new; as to the other results, attention is
paid rather to their interdependence.\bigskip

\noindent
{\bf \S 2. Notation, preliminaries.}
To keep this note halfway self-contained we recall
some definitions and notation (which follow \cite{Pf-boundary}).
$X$ and $Y$ denote Banach spaces which are supposed to be real but in the last section we discuss
some straightforward modifications to settle the complex case;
$X^*$ is the dual of $X$, $\ball(X)$ its closed unit ball.
A bounded set of a Banach space always means a norm bounded set.
$\wstAb{D}$ denotes the $w^*$-closure of a set $D\subset X^*$.
If $K$ is a topological space we endow the set of continuous real valued functions on $K$, $C(K)$, with
the Hausdorff topology $\tau_p$ of pointwise convergence on $K$. If $A\subset X$ and $D\subset X^*$ are sets then $A\eing{D}$
denotes the set of restrictions $\{x\eing{D}\mdE x\in A\}$ which is a subset of $C(D)$ where $D$ carries the
\wst-topology inherited from $X^*$. Similarly, $D\eing{A}$ is defined.
Recall that a subset $A$ of a topological space $T$ is called
relatively countably compact in $T$ if each countable set in $A$ has a cluster point in  $T$
and is called relatively sequentially compact in $T$ if each sequence in $A$ has a convergent subsequence
with a limit in $T$.
Recall further that 
a topological Hausdorff space is called angelic if every relatively countably compact subset $A$ is relatively compact and if the
closure and the sequential closure of such a set $A$ coincide.
It follows from the work of Grothendieck that if $K$ is a compact set then $(C(K),\tau_p)$ is angelic
(\cite[p.\ 36]{Floret} or \cite{Cas-Go} for a generalization).
The topology on $X$ of pointwise convergence on $D$ (where $D\subset X^*$) is denoted by $\tau_D$;
this topology is not necessarily Hausdorff (but it is if $D$ is a norming set for $X$).
Consequently we do not assume that the various compactness notions include the Hausdorff property;
however, this precaution applies only in the statement of Theorem \ref{boundary2.0}.

For a bounded subset $A\subset X$ we define
$$\acoinfty(A)=\{\sum_{n=1}^{\infty}\alpha_n x_n
\mdE \alpha_n\in\R,\; \sum_{n=1}^{\infty}\betr{\alpha_n}=1,\, x_n\in A\}.$$
We say that a set $B\subset\ball(X^*)$ is a boundary for $X$ if for each $x\in X$ there
is $b\in B$ such that $b(x)=\norm{x}$.
A set ${\mathcal N}\subset \ball(X^*)$ is called a norming set for $X$ if $\norm{x}=\sup_{x^*\in{\mathcal N}}\betr{x^*(x)}$
for all $x\in X$. A boundary for $X$ is a norming set for $X$.
Note that if ${\mathcal N}$ is norming for $X$ then its absolute convex hull
$aco({\mathcal N})$ is \wst-dense in $\ball(X^*)$ by Hahn-Banach's theorem,
and $\wstAb{{\mathcal N}}$ contains $extr(\ball(X^*))$ by Milman's converse to Krein-Milman's theorem.
Therefore in Theorem \ref{James2.0} below, there would be no loss of generality if we assumed ${\mathcal N}=extr(\ball(X^*))$ right away.
Our references for unexplained Banach space notions are \cite{JohLin} (for boundaries in particular see \cite{FLP-inf-convex}),
\cite{Die-Seq, FHHMPZ}, for James' theorem \cite{Floret}.\medskip

\noindent
{\bf \S 3. The double limit criterion.}
Let us first comment on Grothendieck's double limit criterion.
As already alluded to in the introduction, the vague idea is that if in the dual $X^*$ something is somehow weakly compact
then in $X$ something else is also somehow weakly compact.
For a less vague formulation see the following lemma.
I thank O. Kalenda and J. Spurn\'{y} for indicating relevant references like \cite{CasMarRaj_2006, FHMZ, Granero_2006}.
\begin{lemma}\label{Groth}
Let $F\subset Y$ and $A\subset Y^*$ be bounded sets. If $A\eing{\wstAb{F}}$ (where $\wstAb{F}\subset Y^{**}$)
is relatively countably compact in $(C(\wstAb{F}),\tau_p)$ then 
$F\eing{\wstAb{aco}(A)}$ is relatively compact and relatively sequentially compact in $(C(\wstAb{aco}(A)),\tau_p)$.\smallskip

If, in addition, $A$ is norming for $Y$ then $F$ is relatively weakly compact.
\end{lemma}
\noindent{\em Proof }
Since $\wstAb{F}$ is \wst-compact, the relative countable compactness of $A\eing{\wstAb{F}}$ and the easy part of
the double limit criterion \cite[p.\ 11]{Floret} or \cite[Cor.\ 2.5]{CasMarRaj_2006} yield that
\bgl
\lim_m \lim_n y_m^*(y_n)=\lim_n \lim_m y_m^*(y_n)
\label{gl_Prop}
\egl
(where all involved limits are supposed to exist)
for all sequences $(y_n)$ in $F$ and $(y_m^*)$ in $A$ (hence also in $-A$).
By \cite[Th.\ 3.3]{CasMarRaj_2006}
(which we read with $Z=[-M,M]$ where $M=\sup_{y\in F, y^*\in\wstAb{aco}(A)}\betr{y^*(y)}$, $K=F$, $\eps=0$, $H=(A\cup (-A))\eing{F}\subset Z^K$),
\eqref{gl_Prop} holds for sequences $(y_m^*)$ in $aco(A)=co(A\cup(-A))$ too.
By \cite[p.\ 12]{Floret} or \cite[Cor.\ 2.5]{CasMarRaj_2006}
(which we read with
$D=aco(A)$, $X=\wstAb{D}$, $H=F\eing{\wstAb{D}}$) it also holds for all sequences $(y_m^*)$ in $\wstAb{aco}(A)$
and the relative compactness of $F\eing{\wstAb{aco}(A)}$ in $(C(\wstAb{aco}(A)),\tau_p)$ obtains.
Relative sequential compactness follows from angelicity of $(C(\wstAb{aco}(A)),\tau_p)$.

If $A$ is norming for $Y$ then $\wstAb{aco}(A)=\ball(Y^*)$ and
$F$ is relatively weakly compact by \wst-compactness of $\wstAb{F}$ and by Banach-Dieudonn\'e's theorem.
\ebew\bigskip

\noindent
{\bf \S 4. Godefroy's conjecture on James' theorem.}\\
Recall the main theorem of \cite{Pf-boundary}:
\begin{theorem}[\textquotedblleft{}boundary2.0{}\textquotedblright]
\label{boundary2.0}
Let $A$ be a bounded set in $X$ and $F\subset X^*$ be a bounded subset
such that for each $x\in\acoinfty(A)$ there is $f\in F$ such that $f(x)=\sup_{g\in \wstAb{co}(F)}g(x)$.\\
If $A$ is $\tau_F$-countably compact then it is $\tau_{\wstAb{co}(F)}$-sequentially compact.
\end{theorem}\noindent
Note that, though stated more generally here,  the just quoted theorem follows
from \cite[Th.~ 9]{Pf-boundary} (with $G=co(F)$) and the remark preceding it. This slightly better statement and Lemma \ref{Groth}
enable us to prove the following a bit technical proposition
by following the proof of James' theorem
in \cite[Cor.~ 10]{Pf-boundary} almost word by word:
\begin{prop}[\textquotedblleft\ James2.1\ \textquotedblright]\label{James2.1}
Let $F\subset Y$ and $A\subset Y^*$ be bounded sets.
If each
\bglst
y^*\in\acoinfty(\wstAb{A})\subset Y^*
\eglst
attains its supremum on $F$ then
$F\eing{\wstAb{aco}(A)}$ is relatively compact and relatively sequentially compact in $(C(\wstAb{aco}(A)),\tau_p)$.
\end{prop}
\noindent
{\em Proof}:
Suppose each $y^*\in\acoinfty(\wstAb{A})$ attains its supremum on $F$.
Let us first note that we may and do assume that $Y$ is the norm closed linear span of $F$
which allows us to identify the (Hausdorff) \wst-topology and $\tau_F$ on bounded sets of $Y^*$.
Set $X=Y^*$. The set $\wstAb{A}$ is \wst-compact in $Y^*$ hence $\tau_F$-compact in $X$.
Theorem \ref{boundary2.0} yields that $\wstAb{A}$ is $\tau_{\wstAb{F}}$-sequentially compact
(where $\wstAb{F}\subset Y^{**}$).
Now the conclusion follows from Lemma \ref{Groth}.\ebew\medskip\\
Immediately we obtain a weaker but less technical formulation
(hence the names 'James2.1' and 'James2.0'):
\begin{theorem}[\textquotedblleft\ James2.0\ \textquotedblright]
\label{James2.0}
Let ${\mathcal N}$ be a norming set for a Banach space $Y$.
For a weakly closed bounded subset $F$ of $Y$ to be weakly compact it is
(necessary and) sufficient that each
\bglst
y^*\in\acoinfty(\wstAb{{\mathcal N}})\subset Y^*
\eglst
attains its supremum on $F$.
\end{theorem}
\noindent
{\em Proof}:
Apply Proposition \ref{James2.1} with $A={\mathcal N}$, $\wstAb{aco}(A)=\ball(Y^*)$.
Then the conclusion follows from the \wst-compactness of $\wstAb{F}$ in $Y^{**}$
and Banach-Dieudonn\'e's theorem.\ebew\bigskip\\
In general Theorem \ref{James2.0} does not generalize James' theorem, for example for strictly convex Banach spaces
we have $\acoinfty(\wstAb{{\mathcal N}})=\ball(Y^*)$.
$C(K)$-spaces form a natural class of Banach spaces to which Theorem \ref{James2.0} might be applied,
in particular for $Y=l^{\infty}$ and ${\mathcal N}=extr(\ball(M(\beta\N)))=\wstAb{\mathcal N}$ we obtain
Godefroy's original conjecture \cite[p.~ 46]{G-Murcia}:
\begin{coro}[\textquotedblleft{}James2.0 for $l^{\infty}${}\textquotedblright]
\label{James2.0for_l_infty}
A weakly closed bounded (not necessarily convex) set $F\subset l^{\infty}$ is weakly compact if each normalized
$\mu\in M(\beta\N)$ ($=(l^{\infty})^*$) with countable support (that is to say of the form \eqref{gl01})
attains its supremum on $F$.
\end{coro}\bigskip

\noindent
{\bf \S 5. To and fro between James' theorem and the boundary problem.}\\
We just saw how Godefroy's conjecture on James' theorem follows from Theorem \ref{boundary2.0}, i.e.\ from the more general form
of the solution of the boundary problem; this explains the first two arrows in the following line.\bigskip
\bgl\label{gl02}\\
\fbox{boundary2.0}
\overset{\parbox{5.5ex}{{\scriptsize double\vspace{-1.5ex}\\\vspace{-1.5ex}limit\\ crit.}}}{\longrightarrow}
\fbox{James2.0}\longrightarrow\fbox{James2.0 for $l^{\infty}$}
\overset{\parbox{5.5ex}{{\scriptsize Gant-\vspace{-1.5ex}\\ macher}}}{\longrightarrow}
\fbox{boundary}\longrightarrow\fbox{James reflexive}\nonumber
\egl\\
Almost conversely, the solution of the boundary problem
follows from Godefroy's conjecture -- cf.\ the third arrow in \eqref{gl02} --
as explained in \cite{G-Murcia}
and reproduced in the proof of Theorem \ref{boundary} below.
The last arrow in \eqref{gl02} (where \textquotedblleft{}James reflexive\textquotedblright\
refers to James reflexivity criterion \cite{James-refl}) is folklore, cf.\ the introduction.\\
The following theorem is a special case of Theorem \ref{boundary2.0} with $F=B$, $\wstAb{co}(B)=\ball(X^*)$
via Eberlein-\v{S}mulyan's theorem whence the names 'boundary' and 'boundary2.0'.
\begin{theorem}[\textquotedblleft{}boundary{}\textquotedblright]
\label{boundary}
Let $B$ be a boundary for a Banach space $X$.
Then for a weakly closed bounded subset $A$ of $X$ to be weakly compact it is
(necessary and) sufficient that $A$ is $\tau_B$-compact.
\end{theorem}
\noindent{\em Proof}:\\
By Eberlein-\v{S}mulyan it is enough to show that every sequence $(x_n)$ in $A$ admits a weakly convergent subsequence
which will be proved by showing that the operator
\bglst
\fkt{S}{l^1}{X}{(\lambda_n)}{\sum \lambda_n x_n}
\eglst
is weakly compact.\\
Since $B$ is a boundary for $X$ the set $co(B)$ is $w^*$-dense in $\ball_{X^*}$ and likewise for the images under $S^*$ in $(l^{\infty})^*$.
By Gantmacher's theorem, $S$ is weakly compact if $S^*$ is and this in turn happens if $C=S^*(co(B))$
is relatively weakly compact in $l^{\infty}$.\\
By Corollary \ref{James2.0for_l_infty} it is enough to show that each $\mu\in (l^{\infty})^*$ of the form
\bgl
\mu=\sum_{n=1}^{\infty}\alpha_n\delta_{s_n}
 \quad\mbox{ where }  \sum_{n=1}^{\infty}\betr{\alpha_n}=1,\, s_n\in\beta\N, \mbox{ $\delta_{s_n}$ Dirac-measure at $s_n$} \label{gl01}
\egl
attains its supremum on $C$.
(With this notation, $(\delta_j)_{j\in\N}$ is the canonical basis of $l^1$.)

For each $n\in\N$ there is an ultrafilter $\mathcal{U}_n$ on $\N$ such that
$\delta_{s_n}=w^*-\lim_{j, \mathcal{U}_n}\delta_j$.
By $\tau_B$-compactness of $A$,
$z_n=\tau_B-\lim_{j, \mathcal{U}_n}x_j$ exists in $X$ for each $n\in\N$.
Set $z=\sum_{n=1}^{\infty}\alpha_n z_n$.
Recalling that $x_j=S\delta_j$ we have $S^{**}\mu(b)=b(z)$ for all $b\in B$ by construction.
By hypothesis there is $b_0\in B$ such that $b_0(z)=\norm{z}$ hence
\bglst
\sup\mu(C)=\sup(S^{**}\mu)(co(B))=\sup(S^{**}\mu)(B)=\sup_{b\in B}b(z)=\norm{z}=b_0(z)
\eglst
i.e. $\mu$ attains its supremum on $C$ at $S^*b_0$.\ebew\medskip\\
Remark.
There is something curious in the last proof because we deal with two boundaries at the same time,
the original one which is encoded in $C$ in $l^{\infty}$ and a second and auxiliary one, the extremal points in the dual of $l^{\infty}$
which appears in Theorem \ref{James2.0} if we write ${\mathcal N}=extr(\ball(X^*))$.
Is there a deeper reason behind this?
\bigskip\bigskip\\
Now we modify \eqref{gl02} by strengthening 'James2.0' so to obtain two equivalent statements:
\bgl
\fbox{boundary2.0}
&\overset{\parbox{5.5ex}{{\scriptsize double\vspace{-1.5ex}\\\vspace{-1.5ex}limit\\ crit.}}}{\longleftrightarrow}
&\fbox{James2.1}\label{gl04}
\egl\\
For \textquotedblleft $\longrightarrow$\textquotedblright of \eqref{gl04} see the proof of Proposition \ref{James2.1}.\\
Proof of \textquotedblleft $\longleftarrow$\textquotedblright of \eqref{gl04}:\\
Let $A\subset X$, $F\subset X^*$ be as in Theorem \ref{boundary2.0},
and suppose that $A$ is $\tau_F$-countably compact.
Consider $x^{**}=\sum\alpha_k x_k^{**}\in\acoinfty(\wstAb{A})\subset Y^*$ with $x_k^{**}\in\wstAb{A}$ and $Y=X^*$.
Each $x_k^{**}$ is a \wst-cluster point of $A$ hence by $\tau_F$-compactness of $A$ there is a
$\tau_F$-cluster point $x_k\in A$ which coincides with $x_k^{**}$ on $F$.
Set $x=\sum\alpha_k x_k$.
Then $x$ and $x^{**}$ coincide on $F$.
Since $x\in\acoinfty(A)$ there is $x^*\in F$ such that $x^{**}(x^*)=x^*(x)=\sup x^{**}(F)$.
By Proposition \ref{James2.1}, $F\eing{\wstAb{aco}(A)}$ is relatively compact
in $(C(\wstAb{aco}(A)),\tau_p)$,
in particular $F\eing{\wstAb{A}}$ is relatively compact in $(C(\wstAb{A}),\tau_p)$,
By Lemma \ref{Groth}, $A\eing{\wstAb{co}(F)}$ is relatively sequentially compact in $(C(\wstAb{co}(F)),\tau_p)$.
Via $\tau_F$-compactness of $A$, this translates into  $\tau_{\wstAb{co}(F)}$-sequential compactness of $A$ hence the conclusion.
\ebew\bigskip

\noindent
{\bf \S 6. Complex scalars.}\\
For all results the complex case follows from the real one by routine arguments.

In the complex case, the proof of Lemma \ref{Groth} works with $H=\bigcup_{\theta\in[0,2\pi]}e^{i\theta}A\eing{F}$
instead of $H=(A\cup (-A))\eing{F}$.

As to Theorem \ref{boundary2.0}, if $X$ is complex then we adapt the statement by substituting
\textquotedblleft$\betr{f(x)}=\sup_{g\in\wstAb{co}(F)}\betr{g(x)}$\textquotedblright\
for \textquotedblleft$f(x)=\sup_{g\in\wstAb{co}(F)}g(x)$\textquotedblright.
The passage from complex to real works like in \cite[Th.~5.10]{FLP-inf-convex}:
The map $x^*\mapsto\Real x^*$ defines an $\R$-linear isometry from $X^*$ onto $(X_{\R})^*$
where $X_{\R}$ is $X$ considered as a real Banach space.
We define $\tilde F=\{\Real(\lambda f)\mdE \betr{\lambda}=1, f\in F\}$ and
$\tilde G=\{\Real(\lambda g)\mdE \betr{\lambda}=1, g\in \wstAb{co}(F)\}$ in $(X_{\R})^*$.
Then $A$ is $\tau_F$-countably compact (respectively $\tau_{\wstAb{co}(F)}$-sequentially compact)
if and only if it is $\tau_{\tilde F}$-countably compact (respectively $\tau_{\tilde G}$-sequentially compact).
The set $\tilde G$ is \wst-closed in $(X_{\R})^*$, and $\sup_{g\in\wstAb{co}(F)}\betr{g(x)}=\sup_{\tilde g\in\tilde G}\tilde g(x)$
for all $x\in X$. For $x\in\acoinfty(A)$ let $f\in F$ be such that $\betr{f(x)}=\sup_{g\in\wstAb{co}(F)}\betr{g(x)}$.
Set $\tilde f=\Real(\frac{\betr{f(x)}}{f(x)}f)$ if $f(x)\neq0$ and $\tilde f=\Real f$ if $f(x)=0$.
By construction, $\tilde f\in\tilde F$ and $\tilde f(x)=\sup_{\tilde g\in\tilde G}\tilde g(x)$.
Now the complex result follows from the real one.

Proposition \ref{James2.1} and Theorem \ref{James2.0} are valid for complex scalars if \textquotedblleft attains its supremum on $F$\textquotedblright\
is replaced by \textquotedblleft admits an $f_0\in F$ such that $\betr{y^*(f_0)}=\sup_{f\in F}\betr{y^*(f)}$\textquotedblright.
Analoguous remarks hold for Corollary \ref{James2.0for_l_infty} and Theorem \ref{boundary}.
(It has already been mentioned in \cite[Th.~ 2.19]{Hardtke} that Theorem \ref{boundary} holds for complex Banach spaces.)

\end{document}